\documentstyle{article}
\title{Canonical Forms of Matrices 
Determining Analytical Manifolds}
\author{Kostadin Tren\v{c}evski and Samet Kera}
\date{}

\begin{document}
\maketitle

\begin{abstract}
In this paper many classes of sets of matrices 
with entries in $F$ $(F\in \{ R,C,H\} )$
are introduced. Each class with the corresponding topology  
determines real analytical, complex or 
symplectic manifold for $F=R$, $F=C$ or $F=H$ respectively.  
Any such family is called to be a set of canonical forms of matrices. 
The construction of such canonical forms of matrices is  
introduced inductively. First basic canonical forms are introduced, 
and then two operations for obtaining new canonical forms 
by using the old canonical forms are introduced. 
All such manifolds have the property that 
each of them can be decomposed into cells of type $F^i$. 

\begin{description}
\item[] Mathematics Subject Classification: 53C25
\item[] Keywords: canonical forms of matrices, analytical 
manifolds, cell decomposition 

\end{description}
\end{abstract}

\section{Introduction}

In the recent paper \cite{T1} are introduced two different 
classes of canonical forms of matrices over a field $F$. 
It easily can be generalized for the quaternions $H$. 
Note that in \cite{T1} and also in this paper the term 
{\em matrix in canonical form} does not mean any reduction of a 
given matrix in a special form, but only means that the corresponding matrix
belongs to a given family or set of matrices. The term {\em canonical form}
comes from the example at the end of this section. According 
to the corresponding topology of the classes of canonical forms of matrices 
(cfm) in this paper are obtained real analytical (for $F=R$), 
complex (for $F=C$) 
and symplectic (for $F=H$) 
manifolds. In all cases are obtained manifolds such that in special 
cases are obtained the Grassmann manifolds. 
In the present paper will be described 
inductively a large class of cfm 
yielding to analytical manifolds. 

Each set of given cfm consists of $n\times m$ matrices with the following 
properties. 

$1^o$. The first property tells about some restrictions concerning the 
matrices in the given canonical form (cf), about the first 
zero coordinates of each vector row. This property depends on the 
choice of the cf. 

$2^o$. The second property is fixed for all cf and it states that any two 
different vector rows are orthogonal.

$3^o$. The third property is also fixed for all cf and states that any vector
row in cf must have norm 1 and the first nonzero coordinate is positive 
real number. 

Thus in order to define a class of cf it is sufficient to specify the 
property $1^o$. Note that alternatively 
the vector rows also can be considered as vectors from 
$RP^{m-1}$, $CP^{m-1}$ or $HP^{m-1}$
and then the property $3^o$ should be omitted. 
Note also that $m\ge n$ according to $2^o$. Indeed, if $m<n$, then  that  set 
of cfm is empty set. 



We finish the introduction by the basic example concerning the 
Grassmann manifolds. 

{\em Example.} Let us consider the set of $n\times m$ matrices $(n\le m)$
such that

$1^o$. If $a_1, \cdots ,a_n \in F^m$ are vector rows, then 
$$0\le t(a_1)< t(a_2)< \cdots < t(a_n)< m,$$
where $t(a_i)$ $(i=1,\cdots ,n)$, denotes the number of 
the first zero coordinates of 
the vector $a_i$. 
This set of matrices together with the fixed properties $2^o$ and $3^o$ 
for cf determines the Grassmann manifold $G_{n,m}(F)$ with the known 
topology, consisting of all $n$-dimensional subspaces of $F^m$ 
generated by the vector rows. 

\section{Construction of different classes of cfm and 
the corresponding analytical manifolds} 

First we determine {\em basic canonical forms} as unit 
$n\times n$ matrices for arbitrary $n$. 

Further we introduce two operations over the cfm, such that the result 
is a new canonical form. After introducing the topology on the new cfm
we obtain a new  
analytical manifold (real, complex and symplectic). 
\smallskip 

i) {\em Inner sum.}

Let $C_1, C_2,\cdots ,C_p$ be $p$ matrices in cf, not necessary in 
cf of the same type. Assume that $C_i$ is an $n_i \times m$ matrix 
$(1\le i\le p)$ and let $C$ be the following $n\times m$ block matrix
$$ C = \left [ \matrix{ C_1\cr C_2\cr \cdot \cr \cdot \cr \cdot \cr
C_p \cr }\right ],$$
where $n=n_1+\cdots +n_p$. If any two different vector rows of $C$ are 
orthogonal, we say that $C$ is in cf called {\em inner sum} of the cfm 
$C_1,\cdots ,C_p$. Thus the inner sum of given $p$ matrices does not 
always exist, while the set of matrices of whole cfm is nonempty for
sufficiently large $m$.    
The term "inner" comes from the orthogonality condition
and we can say only "sum" of the cfm. The new canonical form will 
be denoted symbolically by $C_1+\cdots +C_p$. 
Note that if we neglect the 
orthogonality condition we obtain the Cartesian product of cfm
$C_1,\cdots ,C_p$. 

Now we introduce the topology on the new cfm as follows. 
Let us denote by $\{C_1\}$, ...,$\{C_p\}$ the sets of all matrices 
in given cf, and assume that the corresponding topologies are known. 
Then the set of all matrices in the new cf $\{C_1+...+C_p\}$ 
is a subset of the Cartesian topology space 
$\{C_1\}\times ...\times \{C_p\}$. Thus the topology of  
the new cf $\{C_1+...+C_p\}$ we define to be the relative topology 
with respect to the Cartesian topological space 
$\{C_1\}\times ...\times \{C_p\}$. 

If $\tau :\{ 1,\cdots ,p\} 
\rightarrow \{ 1,\cdots ,p\}$ is a permutation, then obviously 
the cfm $C_1+\cdots +C_p$ and $C_{\tau (1)}+ \cdots + C_{\tau (p)}$ 
determine homeomorphic spaces. 

ii) {\em Spreading of a cfm onto a space $\Sigma$.}

Let $\{ C\}$ be a set of all quadratic $n\times n$ 
matrices from a given canonical form.  
Let $k>n$ and 
$\Sigma $ be a subspace of $F^k$ such that 
$n\le \hbox{dim} \Sigma \le k$. 
Specially, $\Sigma $ can be the total space $F^k$. 
{\em Spreading }
$S_{\Sigma }(C)$ 
over $\Sigma$ is defined to be the set of matrices $CX$ where 
$X$ is $n\times k$ matrix in canonical form of the Grassmann manifold 
$G_{n,k}$ (see example in section 1) such that the row vectors of $X$ 
belong to $\Sigma $, i.e. $X\in G_{n,\Sigma }$. Applying this spreading  
for each matrix $C$ from the given cfm we obtain the total set of matrices of 
the spread canonical form. Notice that  
$C_1X_1=C_2X_2$ implies $C_1=C_2$ and $X_1=X_2$ (see \cite{T1}). 

The topology on the set of matrices in the spread cf we define inductively. 
The topology of the set $\{ C\}$ of quadratic $n\times n$ matrices  
must be the topology of inner sum. So assume that it is obtained as a sum 
$C_1+...+C_p$. Then the topology of the set of the spread cfm is 
defined to be the 
topology of the inner sum of the spread cfm of $C_1$,...,$C_p$ 
on $\Sigma$ separately. Hence it is sufficient to determine the topology   
of each spread cfm $\{ C_i\}$ over $\Sigma$. 
Continuing this process of decreasing the number $n$, we should finally to 
determine the topology of spread cfm of the unit matrix (i.e. basic matrix). 
But, the spreading of the unit matrix is the Grassmann manifold, whose 
topology is well known. 
   
Assuming that $M$ is an analytical manifold, then 
the topology of the spreading $S_{\Sigma }(C)$ is such that 
it is an analytical manifold which is bundle over $G_{n,\Sigma }$ with 
projection $\pi :S_{\Sigma }(C)\rightarrow G_{n,\Sigma }$ defined by: 
$\pi (CX)$ is the subspace of $\Sigma $ generated by the vector rows of 
$CX$, i.e. of $X$, and the fiber is the topology space $M$ induced via 
the set of matrices $\{ C\}$ in initial cf.  
Note that a spreading of a cf can be done only on 
a set of square matrices in cf. Thus 
if $\{ C\}$ is a set of $n\times m$ 
matrices in cf, then we can sum with one or more submatrices in cf (example 
in section 1) determining the Grassmann manifolds in order to obtain 
cf on $m\times m$ matrices. Specially we can sum with $m-n$ vector rows   
or we can sum with a set of 
$(m-n)\times m$ matrices from the Grassmann canonical form. 

Now having in mind the topology of any cfm, we are able to prove that 
each canonical form of $n\times m$ matrices $\{ C\}$ 
yields to an analytical manifold. 
The coordinates of any such analytical manifold will be constructed 
in such a way that for any $n\times n$ nonsingular submatrix will be 
constructed a coordinate neighborhood, like for the standard coordinate 
neighborhoods for the Grassmann manifolds. 

First note that any basic cf of unit $n\times n$ vectors 
determines a 0-dimensional manifold, i.e. a point. 
The first operation must be spreading onto $F^k$ and hence in the first
step we obtain the Grassmann manifolds $G_{n,k}$ which are analytical 
(real, complex and symplectic) manifolds. 

Further, suppose that a new cfm is obtained via an inner sum of matrices 
$C_1,\cdots ,C_p$ in the corresponding cf. If these $p$ cfm determine 
analytical manifolds, then the set of new cfm is also an 
analytical manifold. Indeed, 
let $G_{n_1,\cdots ,n_p,m}(F)$ be the manifold consisting of $(p+1)$-tuples 
$(\Pi _1,\Pi _2,\cdots ,\Pi _p,\Pi _{p+1})$ of orthogonal subspaces of 
$F^m$ where $\hbox{dim} \Pi _i=n_i$ for $i=1,2,\cdots ,p$ and 
$\hbox{dim} \Pi _{p+1}=m-n_1 -\cdots -n_p$. Note that this is a flag manifold of 
included subspaces $V_1=\Pi _1, V_2=\Pi _1+\Pi _2, V_3=\Pi_1 + \Pi_2+
\Pi_3, \cdots $. Thus $G_{n_1,\cdots ,n_p,m}(F)$ is an analytical manifold. 
Let us denote by $M_i$ the analytical manifold which is the fiber of the 
projection of the matrices $\{ C_i \}$ on the Grassmann manifold 
$G_{n_i,m}$. Then the new cf of sum of 
cfm determines a space which is a bundle with base $G_{n_1,\cdots ,n_p,m}(F)$
and fiber $M_1\times \cdots \times M_p$. The projection $\pi $ is given
by 
$$
\pi \Bigl (  \left [ \matrix{ C_1\cr C_2\cr \cdot \cr \cdot \cr \cdot \cr
C_p \cr }\right ]\Bigr ) = 
(\Pi_1,\Pi_2, \cdots ,\Pi_p,\Pi_{p+1}),$$
where $\Pi _1$ is generated by the vector rows of $C_1$, $\Pi_2$ is 
generated by the vector rows of $C_2$, $\cdots$ and $\Pi_{p+1}$ is the 
orthogonal complement of $\Pi_1 + \cdots +\Pi_p$ into $F^m$. Thus this space 
of sum of cfm is also an analytical manifold. 
The coordinate neighborhoods of the new cf can be constructed as follows. 
Let $C'$ be any nonsingular $n\times n$ 
submatrix of $C$. Then there exist $p$ submatrices: 
$C'_1$ submatrix of $C_1$ of order $n_1\times n_1$, 
$\cdots$ , $C'_p$ submatrix of $C_p$ of order $n_p\times n_p$, 
where $n=n_1+\cdots +n_p$, such that 

1. these matrices are submatrices of $C'$, 

2. the columns of $C'_1,\cdots ,C'_p$ are distinct, 

3. the matrices $C'_1,\cdots ,C'_p$ are non-singular, and 

4. by deleting the rows and columns of the submatrices $C'_1, C'_2,\cdots ,
C'_i$, the rest $(n_{i+1}+\cdots +n_p)\times (n_{i+1}+\cdots +n_p)$ 
submatrix of $C'$ is nonsingular $(i=1,2,\cdots ,p-1)$. 

Note that such a choice of submatrices $C'_1, C'_2,\cdots ,C'_p$ is possible 
because by generalization of the Laplace decomposition of the determinants it 
holds 
$$\hbox{det} C' = 
\sum \pm \hbox{det} D'_1 \cdot \hbox{det} D'_2 \cdots \hbox{det} D'_p$$
where $D'_j$ is $n_j\times n_j$ submatrix of $C'$ and submatrix of $C_j$ and
such that the columns of $D'_1,\cdots ,D'_p$ are distinct. 

In the next step we choose the coordinate neighborhoods 
as follows. All $n(m-n)$ elements 
which do not belong to $C'$ may be changed to be close to the corresponding
elements of $C$. The same choice is for the elements 
which are simultaneously in the same row as $C'_i$ and in the same column as 
$C'_j$ for $i>j$, i.e. they may be chosen to be close to the corresponding 
elements of $C'$. Also according to the inductive assumption 
the elements of $C'_1,C'_2,\cdots ,C'_p$ can be 
chosen in the corresponding coordinate neighborhoods which they 
induce respectively on $\{C_1\} ,\{C_2\} ,\cdots ,\{C_p\}$. Finally 
according to the properties 1. - 4. we note that the elements 
which belong simultaneously in the same row as $C'_i$ and in the same column
as $C'_j$ for $i<j$ can uniquely be determined such that the row vectors 
of $C$ are orthogonal. Hence we showed that the chosen matrix $C$ can be 
covered by a coordinate neighborhood of $C_1+\cdots +C_p$. The Jacobi 
matrices for the described coordinate neighborhoods are analytical 
functions, because of the inductive assumptions for 
$C_1,\cdots ,C_p$ and the analytical solutions of linear algebraic systems. 

Next we should show how we can associate a coordinate neighborhood for 
any nonsingular $n\times n$ submatrix $C'$ of spread 
$n\times k$ matrix $S_{\Sigma}(C)$. 
Without loss of generality we can suppose that 
$C$ is a square $n\times n$ matrix in former cf and by inductive assumption  
it can be covered with coordinate neighborhoods with 
analytical elements of the Jacobi matrices. 
Since $\{ S_{\Sigma}(C)\}$ is a bundle with base $G_{n,k}$ 
and according to the standard construction for the coordinates induced by 
any nonsingular $n\times n$ submatrix and the fact that the fiber is an 
analytical manifold by the inductive assumption, we obtain the required 
covering. 

Note that inductively it follows that all the manifolds obtained via this 
method of cf are compact. Finally we can resume the previous results in the 
following two theorems. 

{\bf Theorem 2.1.} {\em The set of all matrices in a chosen 
canonical form with the 
introduced topology is a compact analytical 
(real, complex or symplectic) manifold. }

{\bf Theorem 2.2.} {\em Any manifold determined via a cf of $n\times m$ 
matrices is a bundle over the base manifold $G_{n,m}(F)$ 
or $G_{n_1,n_2,\cdots ,n_p,m}(F)$. }


\section{Some examples} 

In this section will be considered some examples of cfm. 

{\em Example 1.} Assume that $n_1=1,\cdots ,n_p=1$  $(n=p)$ 
and let $s_1,\cdots ,s_n\in \{ 0,\cdots ,m-1\}$ be given numbers.  
Then each of the cf of $1\times m$ matrices 
belongs to the space $FP^{m-1}$. For any $i$, let $\Sigma  _i$  be 
the subspace generated by $e_{s_i+1},e_{s_i+2},\cdots ,e_m$. 
Then the inner sum of cf can be described as a set of $n$ orthogonal 
projective vectors, such that the first starts with at least $s_1$ zeros, 
the second starts 
with at least $s_2$ zeros and so on. The dimension of this manifold is 
$$mn - {n(n+1)\over 2} - s_1 -s_2 -\cdots -s_n.$$
Let us consider the special case $m=n$. Then there exists a permutation 
$\tau $ such that $s_{\tau (i)}<i$, because in opposite case the 
corresponding matrix would not be orthogonal. Hence without loss of 
generality we assume that $s_i<i$, $(i=1,2,\cdots ,n)$. Hence there are 
$n!$ such manifolds. Some of them are homeomorphic. Note that specially, 
if $s_1=s_2=\cdots =s_n=0$, we obtain the full flag manifold. 
\smallskip

{\em Example 2.} Now let us consider the following example. Suppose that 
$m_1,\cdots ,m_p$ are fixed positive integers such that $m=m_1+\cdots +m_p$.
Let 

$\Sigma _1$ be the subspace generated by $e_{1},\cdots ,e_{m_1}$, 

$\Sigma _2$ be the subspace generated by 
$e_{m_1+1},\cdots ,e_{m_1+m_2}$, 

$\Sigma _3$ be the subspace generated by 
$e_{m_1+m_2+1},\cdots ,e_{m_1+m_2+m_3}$, 

$\cdots $ 

$\Sigma _p$ be the subspace generated by 
$e_{m_1+\cdots +m_{p-1}+1}, \cdots ,e_m$. 

Then the cf which is an inner sum has the following form 
$$ \left [ \matrix{C_1 & 0 & 0 & \cdots & 0\cr 
0 & C_2 & 0 & \cdots & 0\cr 
0 & 0 & C_3 & \cdots & 0\cr 
\cdot \cr \cdot \cr \cdot \cr 
0 & 0 & 0 & \cdots &C_p\cr }\right ]$$
as a block matrix of type 
$(n_1+n_2+\cdots +n_p)\times (m_1+m_2+\cdots +m_p)$. Obviously the induced 
analytical manifold is 
the Cartesian product $M_1\times M_2 \times \cdots \times 
M_p$ where $M_i$ is the analytical manifold induced by the $i$-th cf 
of type $C_i$ on $n_i\times m_i$ matrix $(n_i\le m_i)$. 

{\em Example 3.} Let $C$ be an inner sum of $C_1, \cdots  ,C_p$ where 
$C_i$ is an $n_i\times m_i$ matrix, where $n_1+n_2 +\cdots +n_p=n$ and 
$m_1+\cdots +m_p=m=n$. Since $n_i\le m_i$, it must be $n_i=m_i$ for 
$i=1,\cdots ,p$, i.e. 
$$ C = \left [ \matrix{ C_1 & 0 & \cdots & 0 \cr 
0& C_2 & \cdots &0 \cr 
\cdot \cr \cdot \cr \cdot \cr 
0 & 0 & \cdots &C_p\cr } \right ]$$
where $C_i$ is $n_i\times n_i $ matrix in  cf.  Since 
$\hbox{det} C_i\neq 0$ for 
$i=1,2,\cdots ,p$, one can verify that the spread $n\times m'$ matrix $C'$ 
is in cf if and only if $C'$ decomposes into block matrices 
$$C' = \left [\matrix{
C_{11} & C_{12} & \cdots &C_{1p}\cr 
C_{21} & C_{22} & \cdots &C_{2p}\cr 
\cdot \cr 
\cdot \cr 
\cdot \cr 
C_{p1} & C_{p2} & \cdots &C_{pp}\cr }\right ]$$
where $C_{ij}$ is an $n_i\times m'_j$ matrix, where 
$m'_1,\cdots ,m'_p$ are not fixed but $m'_1+\cdots +m'_p=m'$, such that

i) $C_{ij}=0$ for $i>j$, 

ii) rank$C_{ii}=n_i$, 

iii) the row vectors of $C'$ are orthogonal, 

iv) the first non-zero coordinate of each vector is a positive real number. 

In special case if $\{ C_1\} ,\cdots ,\{ C_p\}$ are the cf of full flag 
manifolds, we obtain the manifold described in \cite{T1} sect.3. 

{\em Example 4.} 
Let be given $p$ positive integers $n_{1},n_{2},\cdots ,n_{p}$ and let
$n_{1}+n_{2}+\cdots +n_{p}=n\le m$. 
We consider a set of linearly independent vectors
\par
$$
{\bf a}_{11},{\bf a}_{12},\cdots 
,{\bf a}_{1n_{1}},{\bf a}_{21},{\bf a}_{22},\cdots 
,{\bf a}_{2n_{2}}, \cdots  ,{\bf a}_{p1},{\bf a}_{p2},\cdots 
,{\bf a}_{pm_{p}},
$$
of $F^{m}$ and we denote the matrix with these $n$ 
row-vectors by $A$. The matrix $A$ is in a cf if 
\par
i)\qquad $t_{1}<t_{2}< \cdots <t_{p}$, \quad where 
$t_{i}=\min \{ t({\bf a}_{i1}),t({\bf a}_{i2}),\cdots ,t({\bf a}_{in_{i}})\}$
and $t({\bf a})$
denotes the number of the first zero coordinates of {\bf a},
\par
ii) each two different vectors of these $n$ vectors are orthogonal,
\par
iii) each vector row has norm 1 and the first non-zero 
coordinate is positive real number. 

It is not obvious that this cf belongs to the cf 
of matrices introduced inductively in section 2. 
One can prove that this set of matrices can be considered as 
a spreading of $n\times n$ matrices with the same property 
i), ii) and iii). Thus we should consider the case $m=n$.  
It is clear that the set of $n\times n$ matrices in cf is an inner 
sum of $n_1$ vectors in $F^n$ and set of canonical forms of 
$(n_2+\cdots +n_p)\times n$ matrices 
with parameters $n_2,n_3, \cdots ,n_p$, projected on the space 
generated by the vectors $e_2,e_3, \cdots ,e_n$. Hence by induction of $p$ 
we obtain that the considered cfm is  
included in the family of manifolds obtained in section 2. 

Let us consider the following special case for 
$F=R$, $p=2, n_{1}=2, n_{2}=1$ and $m=3$. Then the 
manifold of canonical vectors consists of the following cells
$$C_{1}=\left [\matrix{
x&*&*\cr y&*&*\cr 0&*&*\cr }\right ], \quad 
C_{2}=\left [\matrix{
1&0&0\cr 0&*&*\cr 0&*&*\cr }\right ], \quad 
C_{3}=\left [\matrix{
0&*&*\cr 1&0&0\cr 0&*&*\cr }\right ], $$
where $x,y>0$. The cell $C_{1}$ is homeomorphic to $R\times S^{1}$ 
because for fixed ratio $\lambda =x/y\in R^{+}$ 
it is homeomorphic to $S^{1}$. The cells $C_{2}$ and $C_{3}$ are 
homeomorphic to $S^{1}$. Thus the Euler characteristic of the manifold 
is $\chi =0$. It can be described such that each point consists of 
two orthogonal lines $p$ and $q$ through the coordinate origin in $R^{3}$ 
such that $q$ lies in the $yz$-plane. 
The third line which is orthogonal to 
$p$ and $q$ is uniquely determined by $p$ and $q$. 
This manifold is homeomorphic to the Klein's bottle. Note that if we 
consider the complex and quaternionic cases, then we obtain complex 
and symplectic manifolds - analogs of the Klein's bottle. 

\section{Decomposition into cells} 

In this section we show the existence of cell decomposition 
which is analogous to the Schubert's cell decomposition of the 
Grassmann manifolds. 

By the construction of cf described in section 2 we obtain a large class of 
compact analytical manifolds, three manifolds for each cfm: 
real, complex and symplectic. 
All these manifolds have the following property. 

{\bf Theorem 4.1.} 
{\em All the manifolds obtained via cfm are such that 
they can be decomposed into disjoint cells of type $F^i$. }

{\em Proof.} Note that the base cf determine 0-dimensional 
manifolds and each of them is a point, i.e. $F^0$. 

First let us prove that if the set of cfm $\{ C\}$ satisfies  the 
property of the theorem 4.1, then $\{ S_{\Sigma}(C)\}$ 
satisfies that 
property also. Since $\{ C\}$ is a set of quadratic matrices, $\{ 
S_{\Sigma}(C)\}$
consists of all matrices of type $CX$ where $X$ is matrix of 
the Grassmann manifolds, and the representation is unique. Hence we obtain
that the cells of $\{ CX\}$ are products of the cells of $\{ C\}$ and the 
cells
of $\{ X\}$. The cells of $\{ C\}$ are of type $F^i$ because of the 
inductive assumption and the cells of the Grassmann manifolds 
$G_{n,\Sigma}$ are also 
of that type and in this case the proof is finished. Indeed, $G_{n,\Sigma}$ 
can be decomposed into ${r\choose n}$ 
cells of type $F^i$, where $r=\hbox{dim} \Sigma$. 

Suppose that the manifolds determined by the cf $\{ C_i\}$ 
satisfy the property 
in theorem 4.1. Then we will show that 
the manifold induced by the sum $C_1+\cdots +C_p$ 
also satisfies that property. This reduces to the special case when 
$C_1,\cdots ,C_p$ are spreadings over the corresponding Grassmann 
manifolds with bases $M_i$ which are quadratic matrices. 
Indeed this manifold is a bundle over the base 
$G_{n_1,n_2,\cdots ,n_p,m}(F)$ and the fiber 
$M_1\times M_2 \times \cdots \times M_p$ and moreover the new 
manifold is equivalent (but not necessary homeomorphic) to the Cartesian 
product 
$$G_{n_1,n_2,\cdots ,n_p,m}(F) \times 
(M_1\times M_2 \times \cdots \times M_p).$$
It follows from the fact that the set of matrices for the spreading of 
quadratic $n_i\times n_i$ matrices is the product (which is unique) of 
matrices of $M_i$ and 
the Grassmann manifold $G_{n_i,m}$. Since $M_i$ satisfies the property in 
theorem 4.1, and the base manifold 
$G_{n_1,n_2,\cdots ,n_p,m}(F)$ can be decomposed into 
$${m!\over n_1! n_2! \cdots n_p! (m-n_1-\cdots -n_p)!}$$
cells of type $F^i$, we obtain that the inner sum 
$C_1+\cdots +C_p$ also satisfies the property of the theorem 4.1. 

This completes the proof of the theorem and moreover 
it gives a method for finding all of the cells. $\Box $

In the paper \cite{TK1} is given a decomposition of the full flag manifold 
$G_n(F)$ into $n!$
cells of type $F^i$. Indeed, the following theorem is proved in \cite{TK1}. 

{\bf Theorem 4.2.} {\em The manifold $G_n(F)$ is a disjoint union of $n!$ 
disjoint cells, such that for each sequence 
$(i_1, i_2, \cdots ,i_{n-1})$, for $0\le i_j \le j$ and $1\le j\le n-1$, 
there exists a cell $C_{i_1,\cdots ,i_{n-1}}$ which is homeomorphic to 
$F^{i_1}\times F^{i_2}\times \cdots \times F^{i_{n-1}}$.}

Note that the theorem 4.1 tells nothing about real manifolds because each 
real manifold can be decomposed 
into cells of type $R^i$. But there are complex 
manifolds which can not be decomposed into disjoint cells of type 
$C^i$. For example if the torus $T=S^1\times S^1$ can be decomposed into 
some cells of type $C^1=R^2$ and $C^0=R^0$, then the Euler characteristic
is a sum of such cells and it is positive number, which is a contradiction. 

The cohomology modules for any manifolds constructed 
via the canonical forms in section 2 can be found easily. 
Indeed, we know the cohomology modules for the manifolds 
$G_{n,m}(F)$ and $G_{n_1,\cdots ,n_p,m}(F)$. 
Using the Leray-Hirsch theorem \cite{BT} we can find step by step
all the cohomology modules for any such manifold. Indeed, according to
the theorem 2.2, we have the following theorem. 

{\bf Theorem 4.3.} {\em Let us denote by $P_t(M)$ the polynomial 
$$\dim H^0(M,R) + t \dim H^1(M,R) + 
t^2 \dim H^2(M,R) + \cdots + t^s\dim H^s(M,R),$$ 
for a manifold $M$, where $s=\dim M$. 

a) If $M$ is a real analytical 
manifold obtained via cfm for $F=R$, then 
$P_t(M)$ is a product of polynomials of types 
$$P_t(G_{p,p+q}) = 
{(1-t)(1-t^2)\cdots (1-t^{p+q})\over 
(1-t)(1-t^2)\cdots (1-t^p)(1-t)(1-t^2)\cdots (1-t^q)} ;$$

b) If $M$ is a complex 
manifold obtained via cfm for $F=C$, then 
$P_t(M)$ is a product of polynomials of types 
$$P_t(G_{p,p+q}) = 
{(1-t^2)(1-t^4)\cdots (1-t^{2(p+q)})\over 
(1-t^2)(1-t^4)\cdots (1-t^{2p})(1-t^2)(1-t^4)\cdots (1-t^{2q})} ;$$

c) If $M$ is a symplectic 
manifold obtained via cfm for $F=H$, then 
$P_t(M)$ is a product of polynomials of types 
$$P_t(G_{p,p+q}) = 
{(1-t^4)(1-t^8)\cdots (1-t^{4(p+q)})\over 
(1-t^4)(1-t^8)\cdots (1-t^{4p})(1-t^4)(1-t^8)\cdots (1-t^{4q})} .$$}

\section{About a further generalization} 

Now will be presented a possible generalization for $F=C$, 
which can be a subject of further research. We give the definition and 
the basic results about symmetric products of manifolds. 

Let $M$ be an arbitrary set and $m$ be a positive integer. In  the
Cartesian product $M^{m}$ we define a relation $\approx $ such that 
$(x_{1},\cdots ,x_{m}) \approx  (y_{1},\cdots ,y_{m})$
iff $y_1,y_2,\ldots ,y_m$ is an arbitrary permutation of 
$x_1,x_2,\ldots ,x_m$. This is equivalence relation and 
the equivalence class  represented  by
$(x_{1},\cdots ,x_{m})$ is denoted by $(x_{1},\cdots ,x_{m})/\approx $  
and  the  quotient space  $M^{m}/\approx $  is called 
{\em symmetric product} of $M$ and is denoted by $M^{(m)}$. 

If $M$ is a topological space, 
then the quotient space $M^{(m)}$ is  also a topological 
space. The space $M^{(m)}$ is introduced quite early \cite{S},  but  mainly
it was studied in \cite{W}. 
The space $(R^{n})^{(m)}$ is a manifold only for
$n=2$ \cite{W}. 
If $n=2$, then $(R^{2})^{(m)}=
{C}^{(m)}$   is
homeomorphic to ${C}^{m}$. Indeed, using that ${C}$ is algebraically closed 
field, it is obvious that the mapping $\varphi :{C}^{(m)}\rightarrow 
C^{m}$ defined by
$$
\varphi ((z_{1},\cdots ,z_{m})/\approx ) = (\sigma _{1},\sigma _{2},\cdots 
,\sigma _{m})
$$
is a bijection, where $\sigma _{i}$ $(1\le i\le m)$ is the $i$-th 
symmetric  function of $z_{1},\cdots ,z_{m}$, i.e.
$$
\sigma _{i}(z_{1},\cdots ,z_{m}) = \sum^{}_{1\le j_{1}<j_{2}<...<j_{i}\le m} 
z_{j_{1}}\cdot z_{j_{2}}\cdots z_{j_{i}}.
$$
The mapping $\varphi $ is also a homeomorphism. 
Moreover, $M^{(m)}$ is a complex manifold if 
$M$ is 1-dimensional complex manifold \cite{T2}. 
For example, if $M$ is a sphere, i.e. the complex 
manifold $CP^1$, then $M^{(m)}$ is the projective complex space 
$CP^m$. Using the permutation products it is easy to see how 
$M^{(m)}=CP^m$ decomposes into disjoint cells 
$C^0,C^1,\cdots ,C^m$. Let $\xi \in M$. Then we define 
$((x_1,\cdots ,x_m)/\approx )\in M_i$ if exactly $i$ of the elements $x_1,
\cdots ,x_m$ are equal to $\xi$. Thus 
$$ M^{(m)} = M_0 \cup  M_1 \cup  \cdots \cup  M_m = 
(M\setminus \{\xi \} )^{(m)}  \cup  
(M\setminus \{\xi \} )^{(m-1)} \cup \cdots \cup 
(M\setminus \{\xi \} )^{(0)}  $$
$$ = C^{(m)} \cup  C^{(m-1)} \cup  \cdots \cup  C^{(0)} = 
C^m \cup  C^{m-1} \cup  \cdots \cup  C^0 .$$
Some recent results about symmetric products of manifolds are obtained 
in \cite{TK2,K1}. This theory about symmetric products has an 
important role in  the  theory  of the topological commutative
vector valued groups \cite{TD1,TD2,T4}. 

We mentioned in the section 1 that the property $3^0$ can be omitted
by assuming that the row vectors of the cfm are projective vectors, i.e. 
they are elements of $CP^{m-1}$. 
Indeed, for any such vector $v=(v_1,\cdots ,v_m)\in CP^{m-1}$ we joint a 
polynomial 
$$P(z) = v_mz^{m-1} + v_{m-1}z^{m-2} + \cdots + v_{1} $$
and hence its complex roots $(z_1,\cdots ,z_{m-1})/\approx $ up to 
a permutation. 
Here $z_1$, $\cdots ,$ $z_{m-1}$ $\in C^*=C\cup \{ \infty \}$ such that 
if $v_m=v_{m-1}=\cdots =v_{m-s+1}=0$ and $v_{m-s}\neq 0$, then exactly $s$
of the roots are equal to $\infty $. 

Now instead of the complex manifold $S^2=C\cup \{ \infty \}$ 
we should consider an 
arbitrary 1-dimensional complex manifold $M$. Then any $1\times m$ canonical 
form of matrices induces the complex manifold $M^{(m-1)}$ and it can be 
considered as a projective space over $M$. The idea for generalization  is   
the following. For any cfm the vector rows should be considered as elements  
of $M^{(m-1)}$. The zero initial values of the vector rows correspond to the 
multiplicity of a chosen point $\xi$ on the chosen 2-dimensional surface.  
If one manage to determine the corresponding orthogonality conditions, 
then he/she will obtain a complex manifold which corresponds to the considered 
cfm and the basic 1-dimensional complex manifold $M$. At this moment we know 
only the projective space over given 2-dimensional real surface, which is the 
symmetric product of the surface, but do not know the Grassmann manifold over 
given 2-dimensional real surface.

\bigskip 

\noindent Institute of Mathematics 

\noindent St. Cyril and Methodius University 

\noindent P.O.Box 162, 1000 Skopje, Macedonia 

\noindent e-mail: kostatre@iunona.pmf.ukim.edu.mk


\noindent e-mail: samet$\_$kera@hootmail.com

\end{document}